\documentclass{article}
\usepackage[utf8]{inputenc}
\usepackage[utf8]{inputenc}
\usepackage{latexsym}
\usepackage{amssymb}
\usepackage{authblk}
\usepackage{amsmath}
\usepackage{amsfonts}
\usepackage{amsthm}
\usepackage{graphicx} % Required for inserting images
\usepackage[english]{babel}
\usepackage{amsmath,amsthm}
\usepackage{amsfonts}
\usepackage{indentfirst}
\usepackage{polynom}
\usepackage{hyperref}
\usepackage{subcaption}
\newcommand{\email}[1]{\texttt{#1}}

\title{\textbf{On One Dimensional Advection - Diffusion Equation with Variable Diffusivity}}

\author{Eeshwar Prasad Poudel\thanks{Institute of Science and Technology, Central Department of Mathematics, Tribhuwan University, Nepal; Tri-Chandra Multiple Campus, Nepal (\email{eeshwarpoudel475@gmail.com})},
Pitambar Acharya\thanks{Department of Applied Mathematics, University of Alabama at Birmingham, AL, USA (\email{pacharya@uab.edu})},
Jeevan Kafle\thanks{Institute of Science and Technology, Central Department of Mathematics, Tribhuwan University, Nepal (\email{jeevankafle@yahoo.com})},
Shreeram Khadka\thanks{Institute of Science and Technology, Central Department of Mathematics, Tribhuwan University, Nepal (\email{shreeramkhadka@gmail.com})}}

\begin{document}
\maketitle
\begin{abstract} 

\noindent In this paper, we address a time-dependent one-dimensional linear advection-diffusion equation with Dirichlet homogeneous boundary conditions. The equation is solved both analytically, using separation of variables, and numerically, employing the finite difference method. The computational output includes three dimensional (3D) plots for solutions, focusing on pollutants such as Ammonia, Carbon monoxide, Carbon dioxide, and Sulphur dioxide. Concentrations, along with their respective diffusivities, are analyzed through  3D plots and actual calculations. To comprehend the diffusivity-concentration relationship for predicting pollutant movement in the air, the domain is divided into two halves. The study explores the behavior of pollutants with higher diffusivity entering regions with lower diffusivity, and vice versa, using 2D and 3D plots. This task is crucial for effective pollution control strategies, and safeguarding the environment and public health.
\end{abstract}


\section{Introduction}
\noindent Advection is defined as a transmission of a matter from one place to another inside a moving fluid \cite{zanudinanalytical}. Diffusion can be described as the changes of molecules from high concentration to low concentration because of the driving force \cite{parker2001encyclopedia}. The general scalar transport equations inclluding the phenomenon advection and diffusion are broadly termed as advection-diffusion equation (ADE) \cite{aswin2015comparative}. In nature, the fluid transport takes place due to combined effect of convection and diffusion \cite{ahmed2017analysis}. ADEs play a pivotal role in applied mathematics and widely used for the simulation of different processes in the areas of environmental sciences, bio sciences, chemical engineering and hydrology \cite{kumar2009analytical}. The one-dimensional ADE is derived on the principle of conservation of mass using Fick’s law \cite{kumar2009analytical}, providing some information about the behaviour of an unknown function: the concentration, and contains its rate of change with respect to location \(x\) and time \(t\) \cite{stockie2011mathematics}. If \(C(x,t)\) is the concentration of pollutants, then the one dimensional advection diffusion equation is \cite{zoppou1997analytical}
\begin{equation}
\frac{\partial C}{\partial t} = D \frac{\partial^2 C}{\partial x^2} - u \frac{\partial C}{\partial x} . \label{eq:advection-diffusion}
\end{equation}
where \(u\) is the velocity with which the fluid flows, also known as the mean advection velocity.
\(D\) is a constant called the diffusion coefficient. 
\noindent The ADE (\ref{eq:advection-diffusion}) describes how the change of concentration \(C\) with respect to position influences its change with respect to time. The solution of ADE provides information about the concentration of pollutants at any time and location \cite{strauss2007partial}. 
\section{Solution of Advection-Diffusion Equation}
\noindent The one-dimensional advection-diffusion equation for a system of length L is given by \cite{zoppou1997analytical}
\begin{equation}
    \frac{\partial C}{\partial t} = D \frac{\partial^2 C}{\partial x^2} - u \frac{\partial C}{\partial x} \quad (0 < x < L, \, 0 < t < \infty). \label{eq:adv-diff-short}
\end{equation}
\noindent The general initial condition and the homogeneous Dirichlet boundary conditions under which the specific solution of the above model is planned to be obtained are \cite{strauss2007partial}
\begin{align}
    C(x, 0) &= f(x), \\
    C(0, t) &= C(L, t) = 0.
\end{align}
\noindent To transform the ADE into a pure diffusion equation, let the nontrivial solution of \eqref{eq:adv-diff-short} be given by
\begin{equation}
    C(x, t) = A(x, t) V(x, t). \label{eq:transformation-C}
\end{equation}
The goal is to find a function $A(x,t)$ such that the advection term in \eqref{eq:adv-diff-short}, i.e., $uC_x$, vanishes.

\noindent Plugging equation \eqref{eq:transformation-C} into equation \eqref{eq:adv-diff-short} leads to
\begin{equation}
    V_t = DV_{xx} + \left[2D\frac{A_x}{A} - u\right] V_x + \left[-u\frac{A_x}{A} + D\frac{A_{xx}}{A} - \frac{A_t}{A}\right] V. \label{eq:collected-equation}
\end{equation}
For the pure diffusion equation, the following conditions should hold
\begin{align}
    2D\frac{A_x}{A} - u &= 0, \label{eq:pure-diff-cond-1} \\
    -u\frac{A_x}{A} + D\frac{A_{xx}}{A} - \frac{A_t}{A} &= 0. \label{eq:pure-diff-cond-2}
\end{align}
Equation \eqref{eq:pure-diff-cond-1} after integration yields
\begin{equation}
    A = K(t) \exp\left(\frac{ux}{2D}\right). \label{eq:A-value}
\end{equation}
Where K(t) is a constant of integration with respect to x and is a function of t only. To determine the value of K(t), from the modified equation \eqref{eq:pure-diff-cond-2}, we obtain
\begin{equation}
    A_t = DA_{xx} - uA_x. \label{eq:A-t-equation}
\end{equation}
Using equation \eqref{eq:A-value}, we have
\begin{equation}
    A_t = K_t \exp\left(\frac{ux}{2D}\right) = \frac{K_t A}{K}. \label{eq:A-t-value}
\end{equation}
% ... (previous code)

\noindent Using these values of $A_t$, $A_x$ and $A_{xx}$; equation \eqref{eq:A-t-equation} becomes
\begin{equation}
    \left(\frac{K_t}{K} - \frac{u^2}{4D} + \frac{u^2}{2D}\right)A = 0. \label{eq:equation-15}
\end{equation}
Assuming that A is a nontrivial solution, we have
\begin{equation}
    \frac{K_t}{K} = -\frac{u^2}{4D}. \label{eq:equation-21}
\end{equation}
Integrating equation \eqref{eq:equation-21} with respect to t yields
\begin{equation}
    K(t) = A_0 \exp\left(-\frac{u^2}{4D}t\right). \label{eq:equation-23}
\end{equation}
where $A_0$ is a pure constant that is independent of both $x$ and $t$.

\noindent Using equation \eqref{eq:equation-23}, equation \eqref{eq:A-value} becomes
\begin{equation}
    A(x, t) = A_0 \exp\left(-\frac{u^2}{4D}t\right) \exp\left(\frac{ux}{2D}\right). \label{eq:equation-24}
\end{equation}
\noindent Substituting the value of $A(x,t)$ from equation \eqref{eq:equation-24} into equation \eqref{eq:collected-equation} and omitting all the details of algebra, we obtain
\begin{equation}
    V_t = DV_{xx}, \label{eq:equation-27}
\end{equation}
which is the transformed pure diffusion equation.\\
The transformed initial condition is 
\begin{equation}
   V(x, 0) = \exp\left(-\frac{ux}{2D}\right) f(x). \label{eq:transf-init-cond}
\end{equation}
The transformed Dirichlet boundary conditions are 
\begin{equation}
     V(0, t) = 0, \label{eq:transf-bound-cond-1} 
\end{equation}
\begin{equation}
    V(L, t) = 0. \label{eq:transf-bound-cond-2}
\end{equation}
\noindent Let the nontrivial solution of equation \eqref{eq:equation-27} be given by
\begin{equation}
    V(x, t) = X(x) T(t). \label{eq:separation-of-variables}
\end{equation}
\noindent where X(x) and T(t) are functions of x only and t only, respectively. Using equation \eqref{eq:separation-of-variables}, equation \eqref{eq:equation-27} becomes
\begin{equation}
    X T_t = D X_{xx} T. \label{eq:separated-eq}
\end{equation}
After rearranging the terms, we obtain
\begin{equation}
    \frac{X_{xx}}{X} = \frac{1}{D} \frac{T_t}{T} = \lambda. \label{eq:separated-eq-lambda}
\end{equation}
where $\lambda$ is a constant. Thus, the transformed diffusion equation \eqref{eq:equation-27} can be split into two ordinary differential equations according to the assumption of equation \eqref{eq:separation-of-variables} as 
\begin{align}
    X_{xx} &= \lambda X, \label{eq:ODE-X} \\
    \frac{T_t}{T} &= D \lambda. \label{eq:ODE-T}
\end{align}
\noindent To obtain the solution of equation (23), the corresponding boundary conditions are
\begin{align}
    V(0, t) &= 0 \Rightarrow X(0) T(t) = 0 \Rightarrow X(0) = 0, \label{eq:boundary-cond-1} \\
    V(L, t) &= 0 \Rightarrow X(L) T(t) = 0 \Rightarrow X(L) = 0. \label{eq:boundary-cond-2}
\end{align}
\noindent Depending upon the value of $\lambda$ (zero, positive, or negative), three different cases arise and the non trivial solution of \label{eq:ODE-X} is possible only in the case $\lambda$ is negative.Thus considering $\lambda = - p^2$ and using the corresponding initial and boundary conditions along with the principle of superposition,the general solution of equation \eqref{eq:separation-of-variables} will be the following linear combination 
\begin{equation}
 V(x, t) = \sum_{n=1}^\infty b_n \sin\left(\frac{n\pi x}{L}\right) \exp\left(-D\left(\frac{n\pi}{L}\right)^2 t\right).
\end{equation}
\noindent where $b_n$ are constants and can be calculated by using the initial condition as
\begin{align}
   V(x, 0) = \sum_{n=1}^\infty b_n \sin\left(\frac{n\pi x}{L}\right) \Rightarrow \exp\left(-\frac{ux}{2D}\right) f(x) = \sum_{n=1}^\infty b_n \sin\left(\frac{n\pi x}{L}\right) 
\end{align}
which is a half range Fourier Sine series and
\begin{align}
 b_n = \frac{2}{L} \int_0^L \exp\left(-\frac{ux}{2D}\right) f(x) \sin\left(\frac{n\pi x}{L}\right) \, dx .  
\end{align}

\noindent Finally, putting together all the results obtained in equation (5), we will obtain the analytical solution to the one-dimensional advection diffusion equation as follows
\begin{align}
    \quad  C(x,t) &= \exp\left(-\frac{u^2 t}{2D}\right) \exp\left(\frac{ux}{2D}\right) \sum_{n=1}^\infty b_n \sin\left(\frac{n\pi x}{L}\right) \exp\left(-D\left(\frac{n\pi}{L}\right)^2 t\right) .
\end{align}
\noindent where $b_n$ can be calculated by using equation (29)\cite{farlow1993partial}.

\section{Comparision}

\noindent Consider an example of one dimensional advection diffusion equation as follows
\begin{align}
&\frac{\partial C}{\partial t} = D \frac{\partial^2 C}{\partial x^2} - u\frac{\partial C}{\partial x} ;  0\leq x \leq 1, t\geq 0 \nonumber \\
&BCs: C(0, t) = C(1, t) = 0; t\geq 0 \\
&IC: C(x,0) = sin\pi x ; 0 \leq x \leq 1. \nonumber
\end{align}

\subsection{Analytical Solution}

\noindent The analytical solution of the above equation is
\begin{align}
C(x,t) = exp\left( -\left(\frac{u^2}{4D}  +  \pi^2 D\right)t\right) sin(\pi x)
   \end{align}

\noindent Under the supposition that $u = 3.6\times 10^{-4} m/hr$ and  $D=3.6 \times 10^{-3} m^2/hr$, the concentration of the pollutants at a distance  $x = 0.6m$  from the point source at time $t = 0.6hr$ is $C(x,t) = exp(-0.03554t). sin(\pi x) = 0.93097$ .

\subsection{Numerical Solution}
\noindent One of the important method to make the study of Numerical Analysis more easier is Finite Difference Method. The main idea in the finite difference method is the derivatives appearing in the differential equation and the boundary conditions are replaced by their finite difference approximations and the resulting linear system of equations are solved by any standard procedure\cite{grewal2014numerical}. 
\noindent Subdividing the spatial interval [0, L] into M + 1 equally spaced sample points $x_m = mh$ , the time interval [0, T ] into N +1 equal time intervals $t_n = nk$, and at each of these space-time points by introducing approximations $C(x_m, t_n) \approx V_m^n$\cite{kafle2020numerical}, the forward time central space scheme (FTCSS) of the one dimensional advection diffusion equation is
\begin{align}
V_m^{n+1} &= V_m^n + \alpha[V_ {m + 1}^n - 2 V_m^n +  V_{m-1}^n] -\beta [V_{m +1}^n - V_m^n]\\
where \quad \alpha &= \frac{hD}{k^2},  \quad \beta = \frac{hu}{k} \nonumber
\end{align}
        Assuming length and the time intervals as $h = 0.2$ and $k = 0.2$ along with  $D = 3.6\times10^{-3} m^2/hr$   and $u = 3.6\times10^{-4} m/hr $ , we have $ \alpha =  0.018$ and  $\beta =  0.00036$.
      
        Also, \begin{align}
         V_0^n &= V_M^n = 0  ,  V_m^0 = sin(\pi mh) \nonumber\\
        \therefore  V_0^0 &= 0 , V_1^0 = 0.5878 , V_2^0 = 0.9511 , V_3^0 = 0.9511 , V_4^0 = 0.5878 , V_5^0 = 0 \nonumber
        \end{align}
    \noindent For $m = 1$ and $n = 0$, from equation (33), we get
    \begin{align}
     V_1^1 &= V_1^0 + \alpha[V_2^0 - 2V_1^0 + V_0^0] - \beta [V_2^0 - V_1^0]\nonumber = 0.58388 
\end{align}
\renewcommand{\arraystretch}{1.7}
\begin{table}[ht]
\caption{Table for Numerical Results} % title of Table
\centering % used for centering table
\begin{tabular}{c c c c c c} % centered columns (6 columns)
\hline\hline %inserts double horizontal lines
$V_0^0=0$ & $V_1^0=0.5878$ & $V_2^0=0.9511$ & $V_3^0=0.9511$ & $V_4^0=0.5878$ & $V_5^0=0$ \\ 
$V_0^1=0$ & $V_1^1=0.58361$ & $V_2^1=0.9445$ & $V_3^1=0.9446$ & $V_4^1=0.5841$ & $V_5^1=0$ \\
$V_0^2=0$ & $V_1^2=0.5794$ & $V_2^2=0.9380$ & $V_3^2=0.9377$ & $V_4^2=0.5802$ & $V_5^2=0$ \\
$V_0^3=0$ & $V_1^3=0.5752$ & $V_2^3=0.9315$ & $V_3^3=0.9313$ & $V_4^3=0.5764$ & $V_5^3=0$\\
$V_0^4=0$ & $V_1^4=0.5711$ & $V_2^4=0.9250$ & $V_3^4=0.9250$ & $V_4^4=0.5726$ & $V_5^4=0$ \\
$V_0^5=0$ & $V_1^5=0.5670$ & $V_2^5=0.9185$ & $V_3^5=0.9187$ & $V_4^5=0.5688$ & $V_5^5=0$\\
\hline %inserts single line
\end{tabular}
\label{table:nonlin} % is used to refer this table in the text
\end{table}

\renewcommand{\arraystretch}{1.7}
\begin{table}[ht]
\caption{Table for Analytical Results} % title of Table
\centering % used for centering table
\resizebox{\textwidth}{!}
{
\begin{tabular}{c c c c c c} % centered columns (6 columns) 
\hline\hline %inserts double horizontal lines
\vspace{0.25cm}
C(0, 0) = 0 & C(0.2, 0) = 0.5878 & C(0.4, 0) = 0.9511 & C(0.6, 0) = 0.9511 & C(0.8, 0) = 0.5878 & C(1, 0) = 0 \\
\vspace{0.25cm}
C(0, 0.2) = 0 & C(0.2, 0.2) = 0.58362 & C(0.4, 0.2) = 0.94432 & C(0.6, 0.2) = 0.94431 & C(0.8, 0.2) = 0.58361 & C(1, 0.2) = 0 \\ 
\vspace{0.25cm}
C(0, 0.4) = 0 & C(0.2, 0.4) = 0.57948 & C(0.4, 0.4) = 0.9376 & C(0.6, 0.4) = 0.9376 & C(0.8, 0.4) = 0.57948 & C(1, 0.4) = 0 \\
\vspace{0.25cm}
C(0, 0.6) = 0 & C(0.2, 0.6) = 0.5753 & C(0.4, 0.6) = 0.93099 & C(0.6, 0.6) = 0.93098 & C(0.8, 0.6) = 0.5753 & C(1, 0.6) = 0 \\
\vspace{0.25cm}
C(0, 0.8) = 0 & C(0.2, 0.8) = 0.5713 & C(0.4, 0.8) = 0.9244 & C(0.6, 0.8) = 0.9243 & C(0.8, 0.8) = 0.5713 & C(1, 0.8) = 0 \\
\vspace{0.25cm}
C(0, 1) = 0 & C(0.2, 1) = 0.5672 & C(0.4, 1) = 0.91785 & C(0.6, 1) = 0.91784 & C(0.8, 1) = 0.56725 & C(1, 1) = 0\\
     
\hline %inserts single line
\end{tabular}
\label{table:nonlin} % is used to refer this table in the text
}
\end{table}

\newpage
\begin{figure}[ht]
\subcaptionbox*{Graph of $C(x,t)$ with varying $x$ and fixed $t$ }[.45\linewidth]{%
    \includegraphics[width=\linewidth]{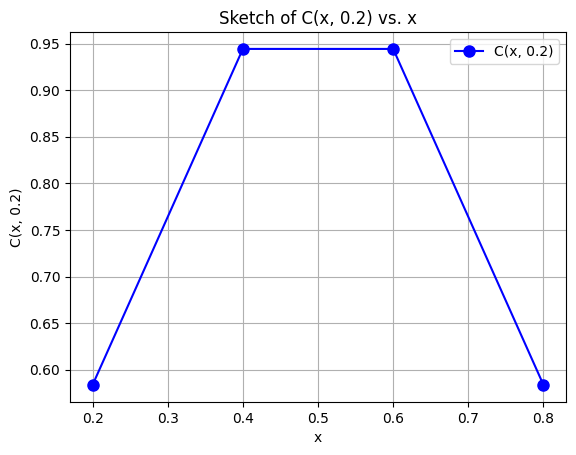}
  }
  \hfill
  \subcaptionbox*{Graph of $C(x,t)$ with time $t$ for various $x$ values}[.45\linewidth]{%
    \includegraphics[width=\linewidth]{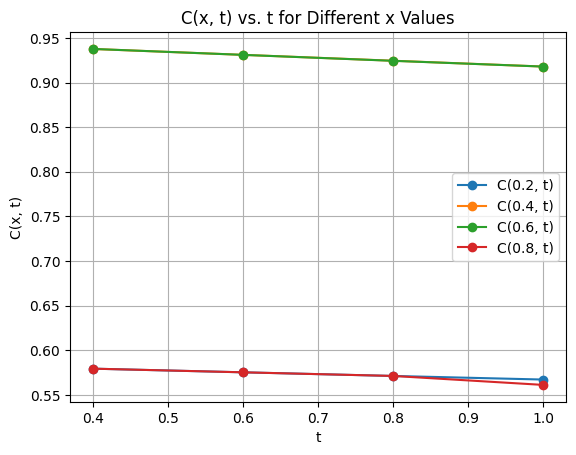}
  }
  
   \caption{Graphical representation}
\end{figure}

\noindent The graph Figure 1(left) illustrates how the value of $C(x, 0.2)$ changes as $x$ varies between $0.2$ and $0.8$. It shows a trend where $C(x,0.2)$ initially increases from $0.58362$ to $0.94432$, remains roughly constant between $x = 0.4$ and $x = 0.6$, and then decreases back to approximately $0.58361$ as $x$ increases further.\\
\noindent The graph Figure 1(right) illustrates how the concentration $C(x,t)$ changes over time for different fixed values of x. It shows how the initial concentrations differ and how they evolve as time progresses.\\
\begin{figure}[ht]
\subcaptionbox*{}[.45\linewidth]{%
    \includegraphics[width=\linewidth]{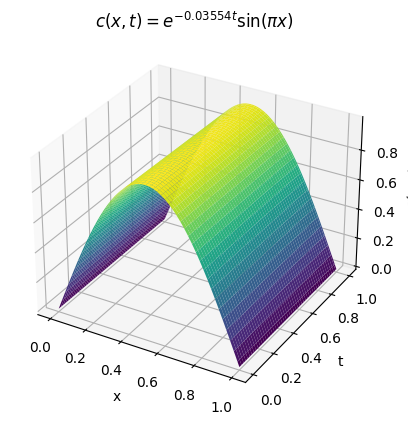}
  }
  \hfill
  \subcaptionbox*{}[.45\linewidth]{%
    \includegraphics[width=\linewidth]{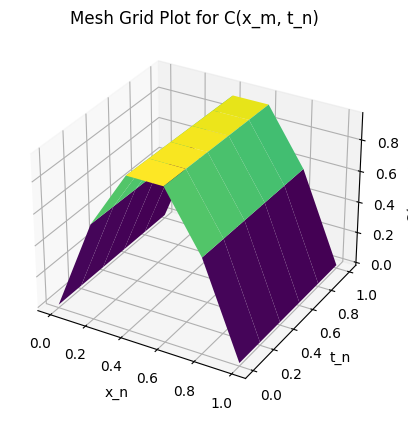}
  }

  \caption{3D plots for Analytical and Numerical solution of ADE}
\end{figure}

\noindent The 3D plot in Figure 2 provides a visual representation of how the function $C(x,t)$ changes with respect to both $x$ and $t$. One can observe the combined effect of exponential decay and oscillatory behavior in the resulting surface. From the mesh grid plot, we see that the function $C(x,t)$ is a product of two terms: The first term, $e^{-0.03554t}$, controls the decay as $t$ increases. As $t$ increases, this term rapidly decreases, indicating that the function becomes smaller as t gets larger. The second term, $\sin(\pi x)$, introduces oscillations along the $x$-axis. The sine wave causes the surface to rise and fall as $x$ varies. As time $t$ increases, the oscillations decrease in magnitude, indicating that the function approaches zero for larger values of $t$. This visual representation helps in understanding the behavior of the function in the $x,t $ space.

\noindent We developed the code for finite difference method. It generated a $3D$ mesh grid plot, (Figure 2) to visualize the behavior of the function $C(x_m,t_n)$ over a $2D$ grid defined by spatial dimension $x_m$ with spatial step size $h=0.2$ and temporal dimension $t_n$  with time step size $k=0.2$. The color of the surface plot at each point $(x_m, t_n)$ represents the value of $C(x_m,t_n)$ at that specific position and time. Brighter (yellow) colors indicate higher values, while darker (blue) colors indicate lower values.\\

\section{Analysis of concentrations of pollutants through 3D plots}

\begin{table}[ht]
\caption{Comparison between concentration of various pollutants} % title of Table
\centering % used for centering table
\resizebox{\textwidth}{!}
{
\begin{tabular}{c c c c c c} % centered columns (6 columns)
\hline\hline %inserts double horizontal lines
Pollutants & Diffusivity $(m^2/hr)$ & Velocity($m/hr$) & $C(x,t)$ & $\alpha$  &$\beta$\\ 
\hline
Ammonia $(NH_3)$ & $7.92\times 10^{-2 }$ & $3.6\times 10^{-4}$ & $e^{-0.78088 t} sin(\pi x)$ & $0.396$ & $0.00036$ \\
Carbon monoxide $(CO)$ & $7.20\times 10^{-2 }$ & $3.6\times 10^{-4} $ & $e^{-0.70989 t} sin(\pi x)$ & $0.360$ & $0.00036$ \\
Carbon dioxide $(CO_2)$  & $5.40\times 10^{-2 }$ & $3.6\times 10^{-4} $  & $e^{-0.53241 t} sin(\pi x)$ & $0.270$ & $0.00036$  \\
Sulphur dioxide $(SO_2)$ & $4.68\times 10^{-2 }$ & $3.6\times 10^{-4}$  & $e^{-0.46142 t} sin(\pi x)$ & $0.234$ & $0.00036$ \\
\hline
\end{tabular}
\label{table:nonlin} % is used to refer this table in the text
}
\end{table}

\begin{figure}[ht]
\subcaptionbox*{Analytical}[.30\linewidth]{%
    \includegraphics[width=\linewidth]{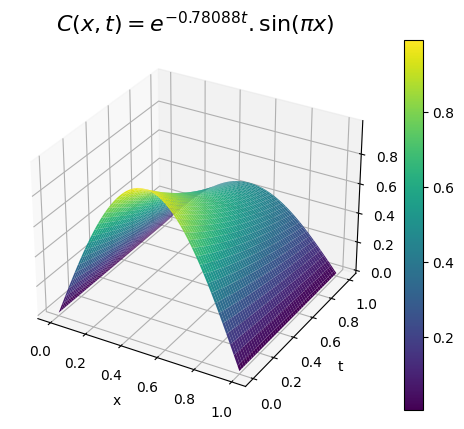}
  }
  \hfill
  \subcaptionbox*{Numerical}[.30\linewidth]{%
    \includegraphics[width=\linewidth]{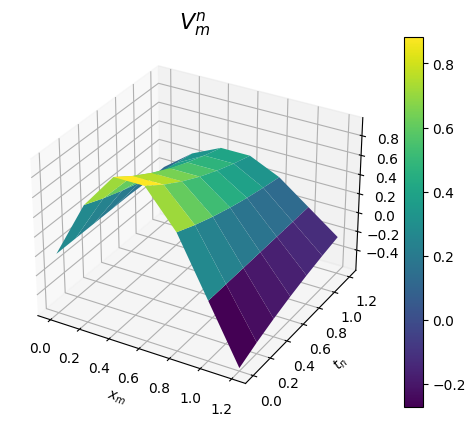}
  }
   \caption{$NH_3 \quad(D = 7.92\times 10^{-2}m^2/hr)$}
\end{figure}

\newpage
\begin{figure}[ht]
\subcaptionbox*{Analytical}[.25\linewidth]{%
    \includegraphics[width=\linewidth]{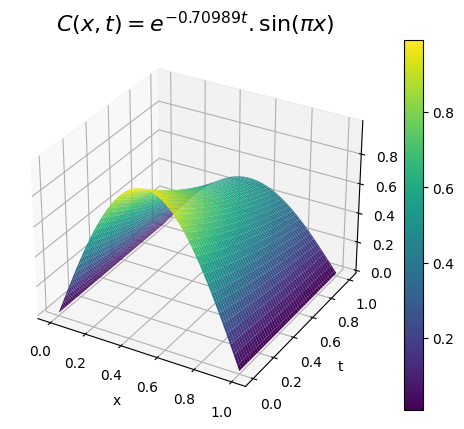}
  }
  \hfill
  \subcaptionbox*{Numerical}[.25\linewidth]{%
    \includegraphics[width=\linewidth]{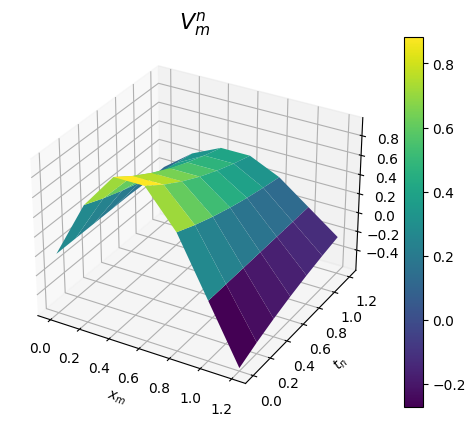}
  }
  \caption{$CO \quad(D = 7.20\times 10^{-2}m^2/hr)$}
\end{figure}

\begin{figure}[ht]
\subcaptionbox*{Analytical}[.25\linewidth]{%
    \includegraphics[width=\linewidth]{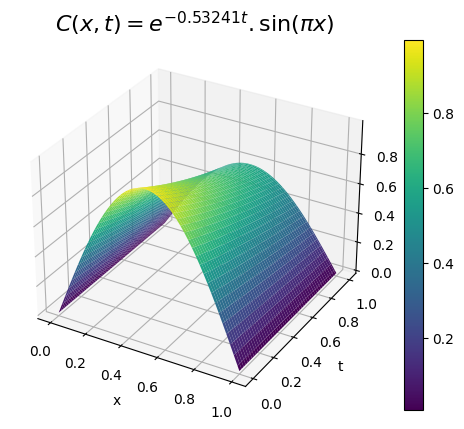}
  }
  \hfill
  \subcaptionbox*{Numerical}[.25\linewidth]{%
    \includegraphics[width=\linewidth]{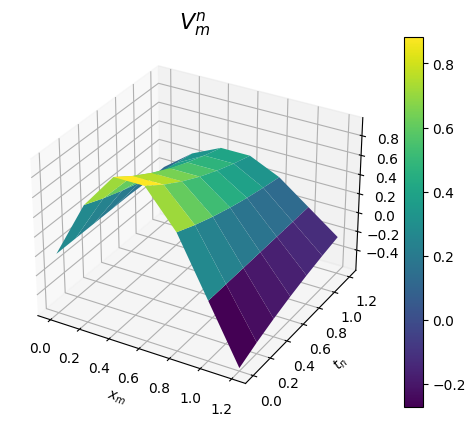}
  }
  \caption{$CO_2 \quad(D = 5.40\times 10^{-2}m^2/hr)$}
\end{figure}

\begin{figure}[ht]
\subcaptionbox*{Analytical}[.25\linewidth]{%
    \includegraphics[width=\linewidth]{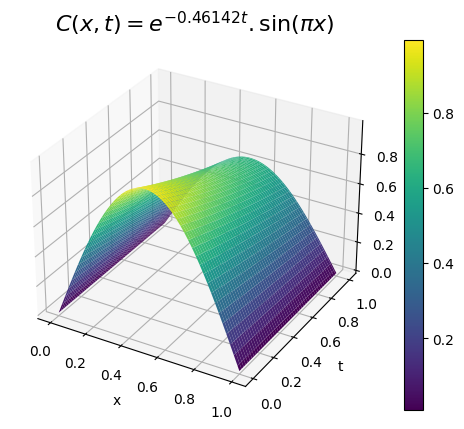}
  }
  \hfill
  \subcaptionbox*{Numerical}[.25\linewidth]{%
    \includegraphics[width=\linewidth]{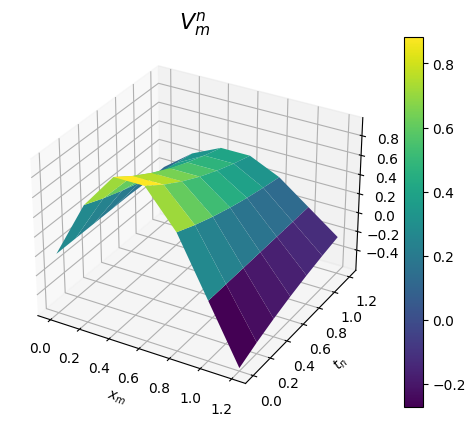}
  }
  \caption{$SO_2 \quad(D = 4.68\times 10^{-2}m^2/hr)$}
\end{figure}

\noindent In this section, we compare the analytical and numerical solution of advection diffusion equation in the case of various pollutants
namely Ammonia - $NH_3 (D = 7.92\times 10^{-2}m^2/hr)$, Carbon monoxide - $CO (D = 7.2\times 10^{-2}m^2/hr)$, Carbondioxide - $CO_2(D = 5.4\times 10^{-2}m^2/hr)$ and Sulphur dioxide - $SO_2 (D = 4.68\times 10^{-2}m^2/hr)$ \cite{reid1977properties, mcketta1964chemical}. From the above four panels of 3D grapical representations, the concentration of pollutants increases as the diffusivity of the pollutants decreases keeping all other factors influencing concentration like molecular size, pressure, nature of the medium etc. constant.

\newpage
\section{Effect of diffusivity of pollutants}

\noindent Dividing the spatial interval [0,1] into two equal halves and assuming the pollutants having some diffusivity enters into the region of other pollutants having different diffusivity, keeping the velocity of both the pollutants constant, the one dimensional advection diffusion equation takes the form as follows: 
\begin{align}
&\frac{\partial C}{\partial t} = D_{1} \frac{\partial^2 C}{\partial x^2} - u\frac{\partial C}{\partial x} ; \quad  0\leq x \leq 0.5, t\geq 0  \\
&\frac{\partial C}{\partial t} = D_{2} \frac{\partial^2 C}{\partial x^2} - u\frac{\partial C}{\partial x} ;  \quad 0.5\leq x \leq 1, t\geq 0 \nonumber 
\end{align}
where $D_1$ and $D_2$ are the diffusivities of pollutants 1 and 2 respectively. Also let's consider:
\begin{align}
&BCs: C(0, t) = C(1, t) = 0; t\geq 0 \\
&IC: C(x,0) = sin\pi x ; 0 \leq x \leq 1. \nonumber
\end{align}

\noindent First, we consider pollutant 1 be Ammonia and pollutant 2 be Sulphur Dioxide. After that we consider pollutant 1 be Sulphur Dioxide ad pollutant 2 be Ammonia. The diffusivity of Ammonia is $7.92\times 10^{-2 } m^2/hr$ and Sulphur Dioxide is $4.68\times 10^{-2 }
m^2/hr$. Referring to figure (7), when the pollutant Ammonia having high diffusivity enters into the region of pollutants Sulphur Dioxide having low diffusivity, the concentration of the pollutants increases rapidly. Also the pollutants having low  diffusivity, like Sulphur Dioxide enters into the region of the pollutants having high diffusivity like Ammonia, the pollutants flows forward slowly and hence the concentration of the pollutants slowly increase due to high diffusivity of right part of the pollutants. 

\begin{figure}[ht]
\subcaptionbox*{(a)}[.45\linewidth]{%
    \includegraphics[width=\linewidth]{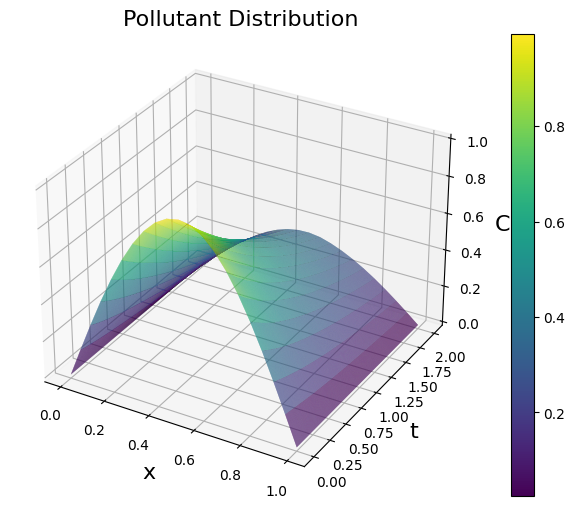}
  }
  \hfill
  \subcaptionbox*{(b)}[.45\linewidth]{%
    \includegraphics[width=\linewidth]{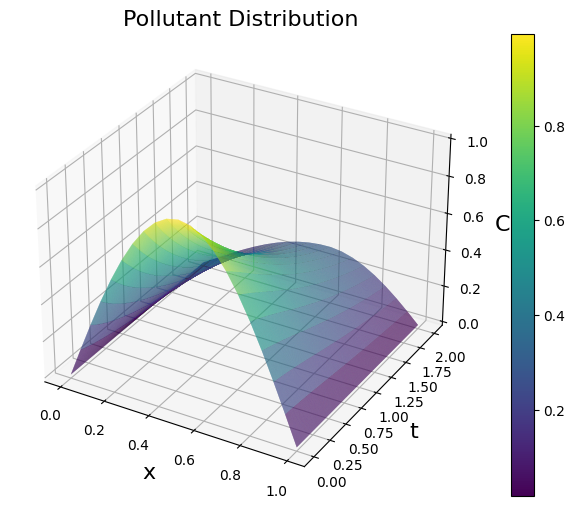}
  }
   \caption{Pollutant distribution with (a) Ammonia at left and Sulphur Dioxide  at right (b) Sulphur Dioxide at left and Ammonia at right}
\end{figure}

\newpage

\begin{figure}[ht]
\subcaptionbox*{(a)}[.55\linewidth]{%
    \includegraphics[width=\linewidth]{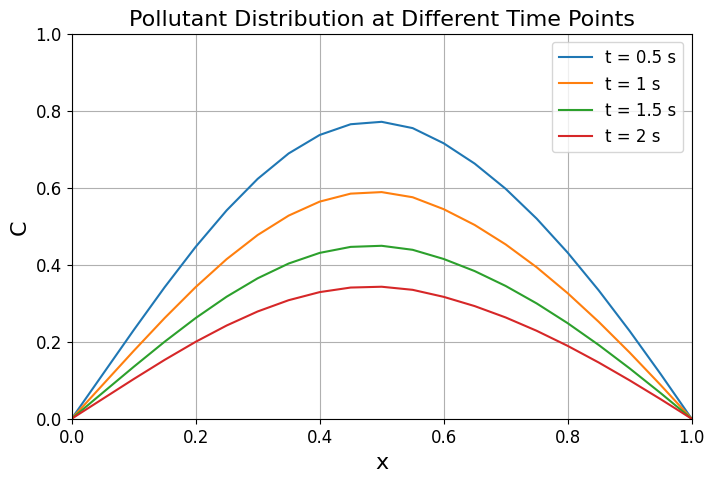}
  }
  \hfill
  \subcaptionbox*{(b)}[.55\linewidth]{%
    \includegraphics[width=\linewidth]{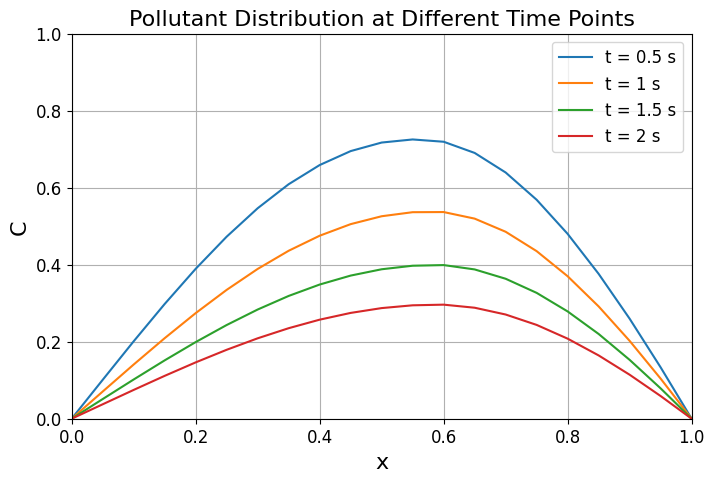}
  }
  \caption{Plot of concentration versus spatial distance at different time levels 0.5s, 1s, 1.5s and 2s with (a) Ammonia at left and Sulphur Dioxide at right (b) Sulphur Dioxide at left and Ammonia at right.}
\end{figure}

\noindent The plot of concentration of pollutants `C' versus `x' at different time levels $t = 0.5s, 1s, 1.5s$, and $2s$ with  (a)Ammonia at left and Sulphur Dioxide at right (b) Sulphur Dioxide at left and Ammonia at right. In each time the maximum concentration `C' are around $x = 0.50m$ and $x= 0.55m$ in figure 8(a) and figure 8(b) respectively. \\

\section{Error Analysis}         
\noindent At  $x = 0.2m$ and $t = 0.2hr$,
 From the table, We have\\ 
Error = $|C_{exact} - C_{approx}|  =0 .00001.\\     
 \noindent So, \% Error=(\frac{0.000010}{0.58362})\times100 = 0.0017\%.$\\
\noindent At  $x = 0.4m$  and $t = 0.6hr$ , From the table , We have \\    
   $\% Error = (\frac{0.000510} {0.93099})\times100 = 0.0547 \%.$ \\
  \noindent  At  $x = 0.6m$ and  $t = 0.6hr$,
 From the table , We have\\    
  \% Error = $(\frac{0.000320}{0.93098})\times100 = 0.0343 \%$.\\
 \noindent At  $x = 0.8m$ and $t = 0.8hr$, From the table, We have \\ 
  \% Error = $(\frac{0.00130}{0.5713}).100 = 0.22 \%$.
  and so on.\\
   \noindent From above calculation, we can observe that the numerical solution by Finite Difference Method for one dimensional advection diffusion equation is very close to its analytical solution and seen to be correct up to four places of decimal. Also, pollutants having higher diffusivity has more error in numerical solution as compared to the analytical solution. Thus, numerical solution is more  appropriate in case of pollutants having lower diffusivity. i.e. FTCS scheme gives better approximation for the pollutants having lower diffusivity as compared to the pollutants having higher diffusivity.

\section {Conclusion} 
\noindent The study presented a one-dimensional advection-diffusion equation as an illustrative example. Analytical solutions were obtained using the separation of variables method, while a numerical solution was derived through the Forward Time Central Space Scheme (FTCSS). Examining various pollutants ($NH_3$, $CO$, $CO_2$ and  $SO_2$ ) with differing diffusivities, it was observed that pollutant concentration increases as diffusivity decreases, assuming constant influencing factors like molecular size, pressure, and medium nature. Treating the atmosphere as a non-homogeneous system, the behavior of pollutants was explored using 2D and 3D plots by equally partitioning high and low diffusivity pollutants in the domain and interchanging them. Notably, when high-diffusivity Ammonia enters the region of low-diffusivity Sulphur Dioxide, pollutant concentration rises rapidly, whereas the entry of low-diffusivity pollutants into high-diffusivity regions results in a gradual increase. Additionally, pollutants with higher diffusivity showed more numerical solution errors compared to analytical solutions, indicating that the FTCS scheme provides a better approximation for pollutants with lower diffusivity.

%\noindent Initially, one dimensional advection diffusion  equation was introduced. Then its analytic solution was obtained by transforming this into pure diffusion equation and then by using separation of variable method with the transformed  initial condition and Dirichlet boundary conditions. We derived Forward Time Central Space Scheme (FTCSS) for one dimensional advection diffusion equation, discussed their stability, consistency and numerical solution was found by using FTCSS. By taking examples of various pollutants like $NH_3$, $CO$, $CO_2$ and  $SO_2$  we studied the analytic and numerical solutions  along with the 3D - plots using python code. We combined the pollutants with high diffusivity and low diffusivity partitioning equally. Further, interchanging the pollutants we studied the Concentration and simulated the 3D and 2D plots and finally calculated errors in several mesh grids.

\clearpage
\bibliographystyle{IEEEtran}
\bibliography{sources.bib}

\end{document}